\documentclass[11pt]{oupau-ppn}
\usepackage{natbib}

\newcommand{\R}{{\mathbb R}}

\newcommand{\Sp}{\mathbb S}
\renewcommand{\H}{\mathrm H}
\renewcommand{\L}{\mathrm L}
\renewcommand{\(}{\left(}
\renewcommand{\)}{\right)}
\newcommand{\be}[1]{\begin{equation}\label{#1}}
\newcommand{\ee}{\end{equation}}
\newcommand{\ird}[1]{\int_{\R^d}{#1}\;dx}
\newcommand{\irdmu}[1]{\int_{\R^d}{#1}\;d\mu_d}
\newcommand{\nr}[2]{\|{#1}\|_{\L^{#2}(\R^d)}}

\newcommand\email[1]{\emph{E-mail:} \textsf{#1}}

\begin{document}
\title{The Euclidean Onofri inequality in higher dimensions}
\shorttitle{The Euclidean Onofri inequality in higher dimensions}
\author{Manuel del Pino\affil{1} and Jean Dolbeault\affil{2}}
\abbrevauthor{del Pino, M., and Dolbeault, J.}
\headabbrevauthor {M. del Pino and J. Dolbeault}
\address{
\affilnum{1} Departamento de Ingenier\'{\i}a Matem\'atica and CMM, UMI CNRS nr.~2807, Universidad de Chile, Casilla 170 Correo~3, Santiago, Chile, \email{delpino@dim.uchile.cl},\\
\affilnum{2} Ceremade, UMR CNRS nr.~7534, Universit\'e Paris-Dauphine, Place de Lattre de Tassigny, 75775 Paris Cedex~16, France, \email{dolbeaul@ceremade.dauphine.fr}.}
\correspdetails{delpino@dim.uchile.cl}

\begin{abstract}\hspace*{-2.8pt} The classical Onofri inequality in the two-dimensional sphere assumes a natural form in the plane when transformed via stereographic projection. We establish an optimal version of a generalization of this inequality in the $d$-dimensional Euclidean space for any $d\ge 2$, by considering the endpoint of a family of optimal Gagliardo-Nirenberg interpolation inequalities. Unlike the two-dimensional case, this extension involves a rather unexpected Sobolev-Orlicz norm, as well as a probability measure no longer related to stereographic projection.\end{abstract}
\maketitle

\medskip\noindent
\emph{Keywords:\/} Sobolev inequality; logarithmic Sobolev inequality; Onofri inequalities; Gagliardo-Nirenberg inequalities; interpolation; extremal functions; optimal constants; stereographic projection\par\smallskip\noindent
\emph{Mathematics Subject Classification (2010):\/} 26D10; 46E35; 58E35
\par\bigskip

\section{Introduction and main result}\label{Sec:Intro}

The \emph{Onofri inequality} as stated in \cite{MR677001} asserts that
\be{Onofri-sphere}
\log\(\int_{\Sp^2}e^v\,d\sigma_2\)-\int_{\Sp^2}v\;d\sigma_2\le\frac 14\,\|\nabla v\|^2_{\L^2(\Sp^2, d\sigma_2)}
\ee
for any function $v\in\H^1(\Sp^2, d\sigma_2)$. Here $d\sigma_2$
 denotes the standard surface measure on the two-dimensional unit sphere $\Sp^2\subset\R^3$, up to a normalization factor $\frac 1{4\,\pi}$ so that $\int_{\Sp^2}1\;d\sigma_2=1$.

\medskip Using stereographic projection from $\Sp^2$ onto $\R^2$, that is defining~$u$~by
\[
u(x)=v(y)\quad\mbox{with}\quad y=(y_1,y_2,y_3)\;,\quad y_1=\frac{2\,x_1}{1+|x|^2}\;,\quad y_2=\frac{2\,x_2}{1+|x|^2}\;,\quad y_3=\frac {1-|x|^2}{1+|x|^2}
\]
for any $x=(x_1,x_2)\in\R^2$, then \eqref{Onofri-sphere} can be reformulated into the \emph{Euclidean Onofri inequality,} namely
\be{Onofri-plane}
\log\(\int_{\R^2}e^u\,d\mu_2\)-\int_{\R^2}u\;d\mu_2\le\frac1{16\,\pi}\,\|\nabla u\|^2_{\L^2(\R^2,dx)}
\ee
for any $u\in\L^1(\R^2,d\mu_2)$ such that $\nabla u\in\L^2(\R^2,dx)$, where
\[
d\mu_2(x):=\frac{dx}{\pi\,(1+|x|^2)^2}
\]
is again a probability measure.

\medskip The purpose of this note is to obtain an (optimal) extension of inequality \eqref{Onofri-plane} to any space dimension. There is a vast literature on Onofri's inequality, and we shall only mention a few works relevant to our main result below. Onofri's inequality with a non-optimal constant was first established by J.~Moser in \cite{MR0301504}, a work prior to that of E.~Onofri,~\cite{MR677001}. For this reason, the inequality is sometimes called the \emph{Moser-Onofri} inequality. We also point out that Onofri's paper is based on an earlier result of T. Aubin,~\cite{MR534672}. We refer the interested reader to \cite{MR2154301} for a recent account on the Moser-Onofri inequality. The inequality has an interesting version in the cylinder $\R\times\Sp^1$, see \cite{MR2437030}, which is however out of the scope of the present work.

\medskip In this note, we will establish that the Euclidean version of Onofri's inequality \eqref{Onofri-plane} can be extended to an arbitrary dimension $d\ge3$ in the following manner. Let us consider the probability measure
\[
d\mu_d(x):=\frac d{|\Sp^{d-1}|}\,\frac{dx}{\(1+|x|^\frac d{d-1}\)^d}\;.
\]
Let us denote
\[
\mathsf R_d(X,Y):=|X+Y|^d-|X|^d-d\,|X|^{d-2}\,X\cdot Y\;,\quad(X,Y)\in\R^d\times\R^d\,,
\]
which is a polynomial if $d$ is even. We define
\[
\mathsf H_d(x,p):=\mathsf R_d\(-\tfrac{d\,|x|^{-\frac{d-2}{d-1}}}{1+|x|^\frac d{d-1}}\,x,\tfrac{d-1}d\,p\)\,,\quad(x,p)\in\R^d\times\R^d\,,
\]
and
\[
\mathcal Q_d[u]:=\frac{\ird{\mathsf H_d(x,\nabla u)}}{\log\(\int_{\R^d}e^u\,d\mu_d\)-\irdmu u}\;.
\]
The following is our main result.
\begin{theorem}\label{Thm:d-Onofri} \emph{With the above notation, for any smooth compactly supported function $u$, we have
\be{Onofri-dEuclidean}
\log\(\int_{\R^d}e^u\,d\mu_d\)-\irdmu u\le\alpha_d\ird{\mathsf H_d(x,\nabla u)}\ .
\ee
The optimal constant $\alpha_d$ is explicit and given by
\[
\alpha_d=\frac{d^{1-d}\,\Gamma(d/2)}{2\,(d-1)\,\pi^{d/2}}\;.
\]
Small multiples of the function
\be{Eigenfunction}
v(x)=\,-\,d\,\frac{x\cdot\mathsf e}{|x|^\frac{d-2}{d-1}\(1+|x|^\frac d{d-1}\)}
\ee
for a unit vector $\mathsf e$ are approximate extremals of \eqref{Onofri-dEuclidean} in the sense that
\[
\lim_{\varepsilon\to 0}\mathcal
Q_d[\varepsilon\,v]=\frac1{\alpha_d}\;.
\]
}\end{theorem}

A rather unexpected feature of inequality \eqref{Onofri-dEuclidean} when compared with Onofri's inequality \eqref{Onofri-plane}, is that it involves an inhomogeneous Sobolev-Orlicz type norm. As we will see below, as a by-product of the proof we obtain a new Poincar\'e inequality in entire space, \eqref{Ineq:Poincare} below, of which the function $v$ defined by \eqref{Eigenfunction} is an extremal.

\begin{example} If $d=2$, $\ird{\mathsf H_2(x,\nabla u)}=\frac 14\int_{\R^2}|\nabla u|^2\,dx$ and we recover Onofri's inequality~\eqref{Onofri-dEuclidean} as in \cite{1101}, with optimal constant $1/\alpha_2=4\,\pi$.
On the other hand, if for instance $d=4$, we find that $H_4(x,\nabla u)$ is a fourth order polynomial in the partial derivatives of $u$, since
$\mathsf R_4(X,Y)=4\,(X\cdot Y)^2+|Y|^2(|Y|^2+4\,X\cdot Y+2\,|X|^2)$.\end{example}

\medskip Extensions of inequality \eqref{Onofri-plane} to higher dimensions were already obtained long ago. Inequality \eqref{Onofri-sphere} was generalized to the $d$-dimensional sphere in \cite{MR1230930,MR1143664}, where natural conformally invariant, non-local generalizations of the Laplacian were used. Those operators are of different nature than the ones in Theorem \ref{Thm:d-Onofri}. Indeed, no clear connection through, for instance, stereographic projection is present. See also \cite{1799948,MR2440885} in which bounded domains are considered.

\medskip Inequality \eqref{Onofri-dEuclidean} determines a natural Sobolev space in which it holds. Indeed, a classical completion argument with respect to a norm corresponding to the integrals defined in both sides of the inequality determines a space on which the inequality still holds. This space can be identified with the set of all functions $u\in\L^1(\R^d,d\mu_d)$ such that the distribution $\nabla u$ is a square integrable function. To avoid technicalities, computations will only be done for smooth, compactly supported functions.

\medskip Our strategy is to consider the Euclidean inequality of Theorem~\ref{Thm:d-Onofri} as the endpoint of a family of optimal interpolation inequalities discovered in \cite{MR1940370} and then extended \hbox{in~\cite{MR1892935}}. These inequalities can be stated as follows.
\begin{theorem}\label{Thm:GN} \emph{Let $p\in(1,d]$, $a>1$ such that $a\leq\frac{p\,(d-1)}{d-p}$ if $p<d$, and $b=p\,\frac{a-1}{p-1}$. There exists a positive constant $\mathsf C_{p,a}$ such that, for any function $f\in\L^a(\R^d,dx)$ with $\nabla f\in\L^p(\R^d,dx)$, we have
\be{Ineq:GN}
\nr fb\leq\mathsf C_{p,a}\,\nr{\nabla f}p^\theta\;\nr fa^{1-\theta}\quad\mbox{with}\;\theta=\tfrac{(a-p)\,d}{(a-1)\,(d\,p-(d-p)\,a)}
\ee
if $a>p$. A similar inequality also holds if $a<p$, namely
\[
\nr fa\leq {\mathsf C_{p,a}}\,\nr{\nabla f}p^\theta\;\nr fb^{1-\theta}\quad\mbox{with}\;\theta=\tfrac{(p-a)\,d}{a\,(d\,(p-a)+p\,(a-1))}\;.
\]
In both cases, equality holds for any function taking the form
\[
f(x)=A\,\(1+B\,|x-x_0|^\frac p{p-1}\)_+^{-\frac{p-1}{a-p}}\quad\forall\,x\in\R^d
\]
for some $(A,\,B,\,x_0)\in\R\times\R\times\R^d$, where $B$ has the sign of $a-p$.}
\end{theorem}
While in~\cite{MR1892935}, only the case $p<d$ was considered, the proof there actually applies to also cover the case~$p=d$, for any $a\in(1,\infty)$.

\medskip For $a=p$, inequality \eqref{Ineq:GN} degenerates into an equality. By substracting it to the inequality, dividing by $a-p$ and taking the limit as $a\to p_+$, we obtain an \emph{optimal Euclidean $\L^p$-Sobolev logarithmic inequality} which goes as follows. Assume that $1<p\le d$. Then for any $u\in\mathrm W^{1,p}(\R^d)$ with $\int_{\R^d}|u|^p\,dx=1$ we have
\[
\int_{\R^d}|u|^p\log|u|^p\,dx\le\frac dp\,\log\left[\beta_{p,d}\int_{\R^d}|\nabla u|^p\,dx\right]\,,\quad\mbox{where}\quad\beta_{p,d}:=\frac pd\,\left(\frac {p-1}e\right)^{p-1}\frac 1{\pi^{\frac p2}}\left[\tfrac{\Gamma\(\tfrac d2+1\)}{\Gamma\(d\,\tfrac{p-1}p+1\)}\right]^{\frac pd}
\]
is the optimal constant. Equality holds if and only if for some $\sigma>0$ and $x_0\in\R^d$
\[
u(x)=\left[\frac 1{2\,\pi^\frac d2}\frac p{p-1}\frac{\Gamma\(\tfrac d2\)}{\Gamma\big(d\,\tfrac{p-1}p\big)}\(\frac p\sigma\)^{d\,\frac{p-1}p}\right]^\frac 1p\;e^{-\frac 1\sigma |x-x_0|^\frac p{p-1}}\quad\forall\; x\in\R^d\,.
\]
This inequality has been established in \cite{MR1957678} when $p<d$ and in general in \cite{MR1990539}; see also \cite{MR2079071,MR2053885}.

\medskip When $p<d$, the endpoint $a=\frac{p\,(d-1)}{d-p}$ corresponds to the usual optimal Sobolev inequality, for which the extremal functions were already known from the celebrated papers by T.~Aubin and G.~Talenti, \cite{MR0448404,MR0463908}. See also~\cite{bliss1930integral,MR0289739} for earlier related computations, which provided the value of some of the best constants.

\medskip When $p=d$, Theorem 1.1 will also be obtained by passing to a limit, namely as $a\to+\infty$. In this way, the $d$-dmensional Onofri inequality corresponds to nothing but a natural extension of the optimal Sobolev's inequality. In dimension $d=2$, with $p=2$, $a=q+1>2$ and $b=2\,q$, it has been recently observed in \cite{1101} that
\[
1\le\lim_{q\to\infty}\mathsf C_{2,q+1}\,\frac{\|\nabla f_q\|_{\L^2(\R^2)}^\frac{q-1}{2\,q}\,\|f_q\|_{\L^{q+1}(\R^2)}^\frac{q+1}{2\,q}}{\|f_q\|_{\L^{2q}(\R^2)}}=\frac{e^{\frac 1{16\,\pi}\,\int_{\R^2}|\nabla u|^2\,dx}}{\int_{\R^2}e^{\,u}\,d\mu_2}
\]
if $f_q=(1+|x|^2)^{-\frac 1{q-1}}\,(1+\frac u{2\,q})$ and $\int_{\R^2}u\;d\mu_2=0$. In that sense, Onofri's inequality in dimension $d=2$ replaces Sobolev's inequality in higher dimensions as an endpoint of the family of Gagliardo-Nirenberg inequalities
\[
\nr f{2q}\le\mathsf C_{2,q+1}\,\nr{\nabla f}2^\theta\,\nr f{q+1}^{1-\theta}
\]
with $\theta=\frac{q-1}q\,\frac d{d+2-q\,(d-2)}$. In dimension $d\ge 3$, we will see below that \eqref{Onofri-dEuclidean} can also be seen as an endpoint of \eqref{Ineq:GN}.

\section{Proof of Theorem~\ref{Thm:d-Onofri}}\label{Sec:Proof}

Assume that $u\in\mathcal D(\R^d)$ is such that $\irdmu u=0$ and let
\[
f_a:=F_a\(1+\tfrac{d-1}{d\,a}\,u\)\;,
\]
where $F_a$ is defined by
\be{Eqn:Optimal}
F_a(x)=\(1+|x|^\frac d{d-1}\)^{-\frac{d-1}{a-d}}\quad\forall\,x\in\R^d\,.
\ee
{}From Theorem~\ref{Thm:GN}, Inequality~\eqref{Ineq:GN}, we know that
\[
1\le\lim_{a\to+\infty}\mathsf C_{d,a}\,\frac{\nr{\nabla f_a}d^\theta\;\nr{f_a}a^{1-\theta}}{\nr{f_a}b}
\]
if $p=d$. Our goal is to identify the right hand side in terms of $u$. We recall that $b=\frac{d\,(a-1)}{d-1}$ and $\theta=\frac{a-d}{d\,(a-1)}$. Using the fact that $F_a$ is an optimal function, we can then rewrite~\eqref{Ineq:GN} with $f=f_a$ as
\[
\frac{\int_{\R^d}|f_a|^\frac{d\,(a-1)}{d-1}\,dx}{\int_{\R^d}|F_a|^\frac{d\,(a-1)}{d-1}\,dx}\le\(\frac{\int_{\R^d}|\nabla f_a|^d\,dx}{\int_{\R^d}|\nabla F_a|^d\,dx}\)^\frac{a-d}{d\,(d-1)}\,\frac{\int_{\R^d}|f_a|^a\,dx}{\int_{\R^d}|F_a|^a\,dx}
\]
and observe that:\medskip

\noindent(i) $\lim_{a\to+\infty}\int_{\R^d}|F_a|^\frac{d\,(a-1)}{d-1}\,dx=\int_{\R^d}\big(1+|x|^\frac d{d-1}\big)^{-d}\,dx=\frac 1d\,|\Sp^{d-1}|$ and
\[
\lim_{a\to+\infty}\int_{\R^d}|f_a|^\frac{d\,(a-1)}{d-1}\,dx=\lim_{a\to+\infty}\int_{\R^d}F_a^\frac{d\,(a-1)}{d-1}\,(1+\tfrac{d-1}{d\,a}\,u)^\frac{d\,(a-1)}{d-1}\,dx=\ird{\frac{e^u}{\(1+|x|^\frac d{d-1}\)^d}}\;,
\]
so that
\[
\lim_{a\to+\infty}\frac{\int_{\R^d}|f_a|^\frac{d\,(a-1)}{d-1}\,dx}{\int_{\R^d}|F_a|^\frac{d\,(a-1)}{d-1}\,dx}=\int_{\R^d}e^{\,u}\,d\mu_d\;.
\]
\\
(ii) As $a\to+\infty$,
\[
\int_{\R^d}|F_a|^a\,dx\approx\frac{2\,a\,\pi^{d/2}}{d^2\,\Gamma(d/2)}\;,\quad \lim_{a\to+\infty}\int_{\R^d}|f_a|^a\,dx=\infty\;,
\]
and
\[
\lim_{a\to+\infty}\frac{\int_{\R^d}|f_a|^a\,dx}{\int_{\R^d}|F_a|^a\,dx}=1\;.
\]
(iii) Finally, as $a\to+\infty$, we also find that
\[
\(\frac{\int_{\R^d}|\nabla f_a|^d\,dx}{\int_{\R^d}|\nabla F_a|^d\,dx}\)^\frac{a-d}{d\,(d-1)}\approx\(1+\tfrac{d\,(d-1)}a\,\alpha_d\ird{\mathsf H_d(x,\nabla u)}\)^\frac{a-d}{d\,(d-1)}\approx\exp\(\alpha_d\ird{\mathsf H_d(x,\nabla u)}\)\;.
\]
Here and above $\ell_1(a)\approx\ell_2(a)$ means that $\lim_{a\to +\infty}\ell_1(a)/\ell_2(a)=1$. Fact (iii) requires some computations which we make explicit next. First of all, we have
\[
\int_{\R^d}|\nabla F_a|^d\,dx=\frac{2\,d^{d-2}\,\pi^{d/2}}{\Gamma(d/2)}\,a^{1-d}\,.
\]
With $X_a:=\(1+\frac{d-1}{d\,a}\,u\)\,\nabla F_a$ and $Y_a:=\frac{d-1}{d\,a}\,F_a\,\nabla u$, we can write, using the definition of $\mathsf R_d$, that
\[
|\nabla f_a|^d=|\nabla F_a|^d\(1+\tfrac{d-1}{d\,a}\,u\)^d+F_a\,|\nabla F_a|^{d-2}\,\nabla F_a\cdot\nabla\(1+\tfrac{d-1}{d\,a}\,u\)^d+\mathsf R_d(X_a,Y_a)\;.
\]
Consider the second term of the right hand side and integrate by parts. A straightforward computation shows that
\[
\int_{\R^d}F_a\,|\nabla F_a|^{d-2}\,\nabla F_a\cdot\nabla\(1+\tfrac{d-1}{d\,a}\,u\)^d\,dx=-\int_{\R^d}|\nabla F_a|^d\(1+\tfrac{d-1}{d\,a}\,u\)^d\,dx-\int_{\R^d}F_a\,\Delta_d F_a\,\(1+\tfrac{d-1}{d\,a}\,u\)^d\,dx
\]
where $\Delta_p F_a=\nabla\cdot(|\nabla F_a|^{p-2}\,\nabla F_a)$ is computed for $p=d$. Collecting terms, we get
\[
\int_{\R^d}|\nabla f_a|^d\,dx=-\int_{\R^d}F_a\,\Delta_d F_a\,\(1+\tfrac{d-1}{d\,a}\,u\)^d\,dx+\ird{\mathsf R_d(X_a,Y_a)}\;.
\]
We may next observe that
\[
a\,\nabla F_a(x)=-\tfrac{d\,a}{a-d}\,|x|^{-\frac{d-2}{d-1}}\,x\(1+|x|^\frac d{d-1}\)_+^{-\frac{a-1}{a-d}}\to -\,d\,\frac{|x|^{-\frac{d-2}{d-1}}\,x}{1+|x|^\frac d{d-1}}\quad\mbox{a.e.}\quad\mbox{as}\quad a\to+\infty\;,
\]
while $a\,\nabla\(1+\tfrac{d-1}{d\,a}\,u\)=\tfrac{d-1}d\,\nabla u$, so that both $X_a=\(1+\frac{d-1}{d\,a}\,u\)\,\nabla F_a$ and $Y_a=\frac{d-1}{d\,a}\,F_a\,\nabla u$ in $\mathsf R_d(X_a,Y_a)$ are of the order of $1/a$. By homogeneity, it follows that
\[
a^d\,\mathsf R_d(X_a,Y_a)\to\mathsf R_d\(-\tfrac{d\,|x|^{-\frac{d-2}{d-1}}}{1+|x|^\frac d{d-1}}\,x,\tfrac{d-1}d\,\nabla u\)=\mathsf H_d(x,\nabla u)\quad\mbox{as}\quad a\to+\infty\;,
\]
by definition of $\mathsf H_d$. Hence we have established the fact that
\begin{multline*}
\int_{\R^d}|\nabla f_a|^d\,dx=-\int_{\R^d}F_a\,\Delta_d F_a\,\(1+\tfrac{d-1}{d\,a}\,u\)^d\,dx+\ird{\mathsf R_d(X_a,Y_a)}\\
=-\int_{\R^d}F_a\,\Delta_d F_a\,\(1+d\,\tfrac{d-1}{d\,a}\,u+o(a^{-1})\)\,dx+a^{-d}\ird{\mathsf H_d(x,\nabla u)}
\end{multline*}
Next we can observe that $-\ird{F_a\,\Delta_d F_a}=\int_{\R^d}|\nabla F_a|^d\,dx$, while $-\lim_{a\to+\infty}a^{d-1}\,F_a\,\Delta_d F_a=d^{d-1}\,|\Sp^{d-1}|\,\mu_d$, so that
\[
-\ird{F_a\,\Delta_d F_a\,u}=a^{1-d}\,d^{d-1}\,|\Sp^{d-1}|\irdmu u+o(a^{1-d})=o(a^{1-d})\quad\mbox{as}\quad a\to+\infty
\]
by the assumption that $\irdmu u=0$. Altogether, this means that
\[
\(\frac{\int_{\R^d}|\nabla f_a|^d\,dx}{\int_{\R^d}|\nabla F_a|^d\,dx}\)^\frac{a-d}{d\,(d-1)}\approx\(1+\frac{\ird{\mathsf H_d(x,\nabla u)}}{a^d\int_{\R^d}|\nabla F_a|^d\,dx}\)^\frac{a-d}{d\,(d-1)}\approx\(1+\tfrac{d\,(d-1)}a\,\alpha_d\ird{\mathsf H_d(x,\nabla u)}\)^\frac{a-d}{d\,(d-1)}
\]
as $a\to+\infty$, which concludes the proof of (iii).

\medskip Before proving the optimality of the constant $\alpha_d$, let us state an intermediate result which of interest in itself. Let us assume that $d\ge 2$ and define $\mathsf Q_d$ as
\[
\mathsf Q_d(X,Y) :=2\,\lim_{\varepsilon\to0}\varepsilon^{-2}\,\mathsf R_d(X, \varepsilon\,Y)=\frac{d^2}{dt^2}\,|X+t\,Y|^d_{\big|t=0}=d\,|X|^{d-4}\left[(d-2)\,(X\cdot Y)^2+|X|^2\,|Y|^2\right]\;.
\]
We also define
\[
\mathsf G_d(x,p):=\mathsf Q_d\(-\tfrac{d\,|x|^{-\frac{d-2}{d-1}}}{1+|x|^\frac d{d-1}}\,x,\tfrac{d-1}d\,p\)\,,\quad(x,p)\in\R^d\times\R^d\,.
\]
\begin{corollary}\label{Cor:Poincare} \emph{With $\alpha_d$ as in Theorem~\ref{Thm:d-Onofri}, we have
\be{Ineq:Poincare}
\int_{\R^d}|v-\overline v|^2\,d\mu_d\le\alpha_d\ird{\mathsf G_d(x,\nabla v)}\quad\mbox{with}\quad\overline v=\irdmu v\;,
\ee
for any $v\in\L^1(\R^d,d\mu_d)$ such that $\nabla v\in\L^2(\R^d,dx)$.}\end{corollary}
\noindent This inequality is a Poincar\'e inequality, which is remarkable. Indeed, if we prove that the optimal constant in~\eqref{Ineq:Poincare} is equal to $\alpha_d$, then $\alpha_d$ is also optimal in Theorem~\ref{Thm:d-Onofri}, Inequality~\eqref{Onofri-dEuclidean}. We will see below that this is the case.

\medskip\noindent\emph{Proof of Corollary~\ref{Cor:Poincare}.} Inequality~\eqref{Ineq:Poincare} is a straightforward consequence of~\eqref{Onofri-dEuclidean}, written with $u$ replaced by $\varepsilon\,v$. In the limit $\varepsilon\to0$, both sides of the inequality are of order $\varepsilon^2$. Details are left to the reader.\hfill \qedsymbol\medskip

To conclude the proof of Theorem~\ref{Thm:d-Onofri}, let us check that there is a nontrivial function $v$ which achieves equality in \eqref{Ineq:Poincare}. Since $F_a$ is optimal for~\eqref{Onofri-dEuclidean}, we can write that
\[
\log\(\ird{|F_a|^\frac{d\,(a-1)}{d-1}}\)=\log\mathsf C_{d,a}+\tfrac{a-d}{d\,(d-1)}\,\log\(\ird{|\nabla F_a|^d}\)+\log\(\ird{|F_a|^a}\)\,.
\]
However, equality also holds true if we replace $F_a$ by $F_{a,\varepsilon}$ with $F_{a,\varepsilon}(x):=F_a(x+\varepsilon\,\mathsf e)$, for an arbitrary given $\mathsf e\in\Sp^{d-1}$, and it is clear that one can differentiate twice with respect to $\varepsilon$ at $\varepsilon=0$. Hence, for any $a>d$, we have
\be{Eqn:Perturbation}
{\textstyle \frac{d\,(a-1)}{d-1}\,\(\frac{d\,(a-1)}{d-1}-1\)}\,\frac{\ird{|F_a|^\frac{d\,(a-1)}{d-1}\,|v_a|^2}}{\ird{|F_a|^\frac{d\,(a-1)}{d-1}}}={\textstyle \frac{a-d}{d\,(d-1)}}\,\frac{\ird{\mathsf Q_d(X_a,\frac{d-1}d\,Y_a)}}{\ird{|\nabla F_a|^d}}+{\textstyle a\,(a-1)}\,\frac{\ird{|F_a|^a\,|v_a|^2}}{\ird{|F_a|^a}}
\ee
with $X_a=\nabla F_a$, $Y_a=\frac d{d-1}\,F_a\,\nabla v_a$ and $v_a:=\mathsf e\cdot\nabla\log F_a$, that is
\[
v_a(x)=-\frac d{a-d}\,\frac{x\cdot\mathsf e}{|x|^\frac{d-2}{d-1}\(1+|x|^\frac d{d-1}\)}\;.
\]
Hence, if $\phi$ is a radial function, we may notice that $\ird{\phi\,v_a}=0$ and
\[
\lim_{a\to+\infty}a^2\ird{\phi\,|v_a|^2}=d^2\ird{\phi(x)\,\frac{|x|^{\frac2{d-1}-2}\,(x\cdot\mathsf e)^2}{\(1+|x|^\frac d{d-1}\)^2}}=d\ird{\phi(x)\,\frac{|x|^\frac2{d-1}}{\(1+|x|^\frac d{d-1}\)^2}}\;.
\]
Since ${\ird{|F_a|^\frac{d\,(a-1)}{d-1}}}=o\({\ird{|F_a|^a}}\)$, the last term in \eqref{Eqn:Perturbation} is negligible compared to the other ones. Passing to the limit as $a\to+\infty$, with $v:=\lim_{a\to+\infty}a\,v_a$, we find that $v$ is given by \eqref{Eigenfunction} and
\[
\(\tfrac d{d-1}\)^2\irdmu{|v|^2}=\alpha_d\ird{\mathsf Q_d\Big(-\tfrac{d\,|x|^{-\frac{d-2}{d-1}}}{1+|x|^\frac d{d-1}}\,x,\tfrac{d-1}d\,Y\Big)}
\]
where $Y:=\frac{d-1}d\,\nabla v$ and where we have used the fact that
\[
d\,(d-1)\,\alpha_d\,\lim_{a\to+\infty}a^d\ird{|\nabla F_a|^d}=1\;.
\]
Since the function $\mathsf Q_d$ is quadratic, we obtain that
\[
\big(\tfrac d{d-1}\big)^2\irdmu{|v|^2}=\alpha_d\ird{\mathsf G_d(x,\tfrac d{d-1}\,\nabla v)}=\alpha_d\,\big(\tfrac d{d-1}\big)^2\ird {\mathsf G_d(x,\nabla v)}\;,
\]
which corresponds precisely to equality in \eqref{Ineq:Poincare} since $v$ given by \eqref{Eigenfunction} is such that $\overline v=0$.

\medskip Equality in~\eqref{Onofri-dEuclidean} is achieved by constants. The optimality of $\alpha_d$ amounts to establish that in the inequality
\[
\mathcal Q_d[u]\ge\frac1{\alpha_d}\;,
\]
equality can be achieved along a minimizing sequence. Notice that
\[
\mathcal Q_d[u]=\frac{\ird{\mathsf H_d(x,\nabla u)}}{\log\(\int_{\R^d}e^u\,d\mu_d\)}\quad\mbox{if}\quad \irdmu u=0\;,
\]
The reader is invited to check that $\lim_{\varepsilon\to0}\mathcal Q_d[\varepsilon\,v]=\frac1{\alpha_d}$. In dimension $d=2$, $v$ is an eigenfunction associated to the eigenvalue problem: $-\Delta\,v=\lambda_1v\mu_2$, corresponding to the lowest positive eigenvalue, $\lambda_1$. The generalization to higher dimensions is given by \eqref{Eigenfunction}. Notice that the function $v$ is an eigenfunction of the linear form associated to $\mathsf G_d$, in the space $L^2(\R^d,d\mu_d)$. This concludes the proof of Theorem~\ref{Thm:d-Onofri}.

\medskip Whether there are non-trivial optimal functions, that is, whether there exists a non-constant function $u$ such that $\mathcal Q_d[u]=\frac1{\alpha_d}$, is an open question. At least the proof of Theorem~\ref{Thm:d-Onofri} shows that there is a loss of compactness in the sense that the limit of $\varepsilon\,v$, \emph{i.e.}~$0$, is not an admissible function for $\mathcal Q_d$.

\medskip\noindent{\small{\bf Acknowlegments.} J.D.~has been supported by the projects \emph{CBDif} and \emph{EVOL} of the French National Research Agency (ANR). M.D. has been supported by grants Fondecyt 1110181 and Fondo Basal CMM. Both authors are participating to the \emph{MathAmSud} network \emph{NAPDE.}}

\bibliographystyle{plainnat}\bibliography{References}\end{document}